\documentclass[12 pt]{article}
\usepackage{amsfonts}
\usepackage[pctex32]{graphics}
\title{The Fair and Random Maximal Division of ``Pizza"}
\author{Floyd E.~Brown\\
Department of Mathematics\\
Science Hill High School
\and
Anant P.~Godbole\\
Department of Mathematics\\
East Tennessee State University}
\begin{document}
\def\qed{\vbox{\hrule\hbox{\vrule\kern3pt\vbox{\kern6pt}\kern3pt\vrule}\hrule}}
\def\ep{\varepsilon}
\def\lr{\left(}
\def\lf{\lfloor}
\def\rf{\rfloor}
\def\lc{\left\{}
\def\rc{\right\}}
\def\rr{\right)}
\def\p{\mathbb P}
\def\v{\mathbb V}
\def\e{\mathbb E}
\def\l{\mathbb L}
\def\lg{{\rm lg}}
\newtheorem{thm}{Theorem}
\newtheorem{lemma}[thm]{Lemma}
\newtheorem{prop}[thm]{Proposition}
\maketitle

\medskip \noindent {\it Keywords:  Minimum standard deviation, minimum mean absolute deviation, central limit theorem, Stein's method.}

\section{Introduction}  The second-named author of this paper has run REU (``Research Experiences for Undergraduates") sites during each summer since 1991.  The project titles have been {\it The Probability Theory of Runs and Patterns; Discrete Probability and Associated Limit Theorems; Probabilistic Methods in Graph Theory, Combinatorics and Number Theory} and {\it Discrete Random Structures}, which will give readers an idea of the {\it kind} of research that is conducted.  During the summers of 1999, 2001 and 2003, the program was enriched by the addition of an RET (``Research Experiences for Teachers") supplement that allowed one, three, and two teachers respectively to work alongside some of the brightest mathematics and statistics majors in the country.  The work described in this paper was conducted by first-named author of this paper, a high school teacher in Tennessee, during the summer of 2001 -- and under the mentorship of the project director.  The need to actively involve K-12 teachers in the process of 
mathematical discovery is now beyond debate.  The increased level of NSF funding for the RET program, and the views of authorities such as Al Cuoco, Director of the Center for Mathematics Education at the Education Development Center, certainly bear testimony to this fact.  We can do no better than to echo Dr. Cuoco's words \cite{cuoco}: 

``There are very few absolutes in education, but there's one thing of which I am absolutely certain:  The best high school teachers are those who have a research-like experience in mathematics.....working for an extended period of time on a hard problem that has no apparent approach or solution has profound effects on how one perceives the nature of the enterprise.  Teachers who have done this type of research are much less likely to think of mathematics as an established body of facts than are teachers who have simply taken a set of courses.  They are more likely to stay engaged in teaching after they start teaching….. And they are much more likely to organize their classes around large investigations rather than low-level exercises....  An ideal teacher preparation program   combines the kind of orchestrated assimilation of the main results in mathematics....with the much messier unstructured explorations that come from working with a mentor and grappling with a research project."
       
 Consider a circular pizza.  If one makes $n$ straight line cuts in the ``standard" way, we end up with $2n$ slices; the answer is, of course, different if one makes the horizontal and vertical cuts favored by some pizza parlors -- but the number of slices is still linear in $n$.  What, however, if one seeks, at the expense of the regularity of the shape of the slices, to {\it maximize} the number of slices?  Key to the solution to this problem is to ensure sequentially that every cut intersects each of the previous ones.  This generates a total of $n^2/2+n/2+1=1+n+{n\choose 2}$ slices with $n$ cuts (\cite{steiner}, \cite{buck}, \cite{kerr}, \cite{wetzel}, \cite{zimmerman}); see Figure 1 below for the case $n=3$.  The maximal number of slices is thus quadratic in 
$n$.  It has been our experience that even elementary school children enjoy a classroom activity based on this fact, and particularly the challenge associated with the successful drawing of a correct solution for $n$ even as small as 5. 

\begin{center}
\vspace{2in}
\centerline{\hbox to 4in{\special{bmp:c:pizza1.bmp x=4in y=2in}}}
\centerline{Figure 1}\
\end{center}

Both aspects of our investigation into this problem are statistical in nature.  In Section 2, we consider the simplest case of $n=3$ and ask a very simple question.  The nature of our process guarantees that the maximal solution will {\it not} lead to slices of equal size.  How might one make the cuts, however, so that the slices are of as equal size as possible?  In other words, if we denote the areas of the slices of a unit pizza by $A_1,\ldots,A_7$, how should the cuts be made so as to minimize $\v\{A_1,\ldots,A_7\}$?  The surprising answer is revealed in what we affectionately labeled the ``Every Seventh Child Starves" principle, which dictates that the variance is minimized in a class of solutions when three ``regular", i.e., radial cuts are made, leading to six slices of area $\pi/6$ and one of size 0.  This is clearly unacceptable from an ethical viewpoint and we next investigate which other measures might be used instead so as to give a more acceptable solution.  

In Section 3, we investigate what happens if $n$ cuts are made so as to potentially maximize the number of slices, but with a blunt knife -- so that the actual number $R=R_n$ of slices is a random variable whose behavior is determined by $p=p_n$, the probability of any cut being successful.  We compute $\e(R)$ and $\v(R)$ and use Stein's method of normal approximation to prove that a normalized version of $R$ is close to normal in distribution as $n\to\infty$.

The research in this paper is characterized by three important pedagogical facts:  First, the problem emanates as an attractive pizza-cutting problem that we have successfully used as an ``emerging pattern problem" as early as at the second grade level.  Second, the maxim that ``minimum possible variance" is equivalent to the data points being ``as equal as possible" is used, together with standard high-school AP-level calculus, to solve the corresponding ``fair division" problem with three cuts.  We have presented this solution to students in our Mathematical Statistics and Calculus classes.  Also, graph paper and rulers have been used to empirically guess the answer to the minimum SD problem at an even lower level -- with each student in a large {\it elementary statistics} class drawing the appropriate picture haphazardly, thus creating his/her own data set with seven numbers, and finding the resulting SD -- and with the results of the class being pooled.  Last but not least, routine calculations at a mathematical-statistics level are used in this paper, in conjunction with a ready-to-use version of Stein's method of normal approximation for sums of dependent random variables, to solve the random division problem with $n$ cuts, $n\to\infty$.  Portions of this material have been used to create non-standard classroom examples, of the computation of moments, for upper level students.

The material in this paper is similar in spirit, but quite different in content, from that in \cite{benfield}

\section{Fair Division:  When is Standard Deviation Inappropriate?}  Figure 1 illustrates the general configuration we believe will direct us to a solution.  We wish to determine a function for the area of each region as a function of the length of arc 
$AB$, denoted by $x$.  A few assumptions first:\\ 
\indent 1) The triangle in the center is equilateral with its center coinciding with the center of the circle, denoted by $O$.  As a result, the chords are congruent and an obvious symmetry emerges with regions $ABE \cong HDF\cong CGI$ and $AEFH\cong BEGI\cong DFGC$;\\
\indent 2) For convenience the circle has a radius of one unit.  For discussion purposes we will refer to the three types of regions as the triangle, the ``circular triangles", and the ``circular trapezoids";\\
\indent 3) $x$ lies in the interval $[0, {\pi\over 3}$].\\

Given our assumptions, let's determine a function for the area of each of the three types of regions as functions of $x$, the length of arc $AB$.  This calculation involves simple   trigonometry but is surprisingly and unexpectedly non-trivial.  We have

$$
m\angle FEG=m\angle AEB={\pi\over3};$$
and
$$m\angle BEO={{5\pi}\over6};\ m\angle EOB={x\over2};\ {\rm and}\ m\angle EBO={{\pi\over6}-{x\over2}}.$$
Also, the law of sines yields, since $OB=1$,
$$
{1\over {\sin{{5\pi}\over6}}}={EB\over {\sin{x\over2}}},\ {\rm thus}\ EB=2\sin{x\over2},$$
and
$${1\over {\sin{{5\pi}\over6}}}={EO\over {\sin{\lr{\pi\over6}-{x\over2}\rr}}}\ {\rm thus}\  EO=2\sin\lr{{\pi\over6}-{x\over2}}\rr.$$
It follows that
\begin{eqnarray*}
\triangle BEO&=&{1\over2}\lr{2\sin{x\over2}}\rr \lr{2\sin{\lr{{\pi\over6}-{x\over2}}\rr}}\rr\sin{{5\pi}\over6}\\
&=&\sin\lr{x\over2}\rr\sin\lr{{\pi\over6}-{x\over2}}\rr.\\
\end{eqnarray*}
Also , the area of sector $AOB$ equals
${x\over2}$ and thus the area
of a circular triangle equals
${x\over2}-2\sin\lr{x\over2}\rr\sin\lr{\pi\over6}-{x\over2}\rr.$
Next, the area of the triangle is computed by noting that its altitude equals
$3\sin\lr{\pi\over6}-{x\over2}\rr,$ and that ${1\over2}FG={3\over{\sqrt {3}}}\sin{\lr{\pi\over6}-{x\over2}\rr}.$
It follows that the area of the triangle equals ${3 {\sqrt {3}}}\sin{^{2}}{\lr{\pi\over6}-{x\over2}\rr}.$  Finally, the
area of circular trapezoid is determined by subtraction to equal
\begin{eqnarray*}
&&{}{1\over3}\lc\pi-3\Big[{x\over2}-2\sin\lr{x\over2}\rr\sin\lr{\pi\over6}-{x\over2}\rr\Big]-3\sqrt{3}\sin^{2}\lr{\pi\over6}-{x\over2}\rr\rc\\
&=&{{\pi}\over3}-{x\over2}+2\sin\lr{x\over2}\rr\sin\lr{\pi\over6}-{x\over2}\rr-\sqrt{3}\sin^{2}\lr{\pi\over6}-{x\over2}\rr.
\end{eqnarray*}
Our three area functions will now be defined as follows: Area of the triangle $=\alpha_1$, area of each circular triangle $=\alpha_2$, and area of each circular trapezoid $=\alpha_3$:

\begin{eqnarray*}
\alpha_1(x)&=&{3 {\sqrt {3}}}\sin{^{2}}{\lr{\pi\over6}-{x\over2}\rr}\\
\alpha_2(x)&=&{x\over2}-2\sin\lr{x\over2}\rr\sin\lr{\pi\over6}-{x\over2}\rr\\
\alpha_3(x)&=&{{\pi}\over3}-{x\over2}+2\sin\lr{x\over2}\rr\sin\lr{\pi\over6}-{x\over2}\rr-\sqrt{3}\sin^{2}\lr{\pi\over6}-{x\over2}\rr.
\end{eqnarray*}

\medskip

\noindent
{\bf ANALYSIS OF THE PROBLEM}
[Note: The calculations and graphs from this point were done using the MAPLE software package.  Much of the ``standard" calculus involved {\it can} be done by hand, but is rather tedious.]

Since $(A_1+A_2+\ldots+A_7)/7=\pi/7$, the standard deviation function is easily seen, on using the identity $\sin u=\cos(\pi/2-u)$, to be given by
\begin{eqnarray*}
\sigma(x)&=&{1\over7}\sqrt \Big[21{\lr{x\over2}-2\sin{\lr{x\over2}\rr}\cos{\lr{\pi\over3}+{x\over2}\rr}\rr}^2+189\cos^4{\lr{\pi\over3}+{x\over2}\rr}\\
&&{}+21\lr{\pi\over3}-{x\over2}+2\sin\lr{x\over2}\rr \cos{\lr{\pi\over3}+{x\over2}\rr}-\sqrt{3}\cos^{2}{\lr{\pi\over3}+{x\over2}\rr}\rr^2-\pi^2 \Big];
\end{eqnarray*}
\bigskip
\hfil\scalebox{1}{\includegraphics{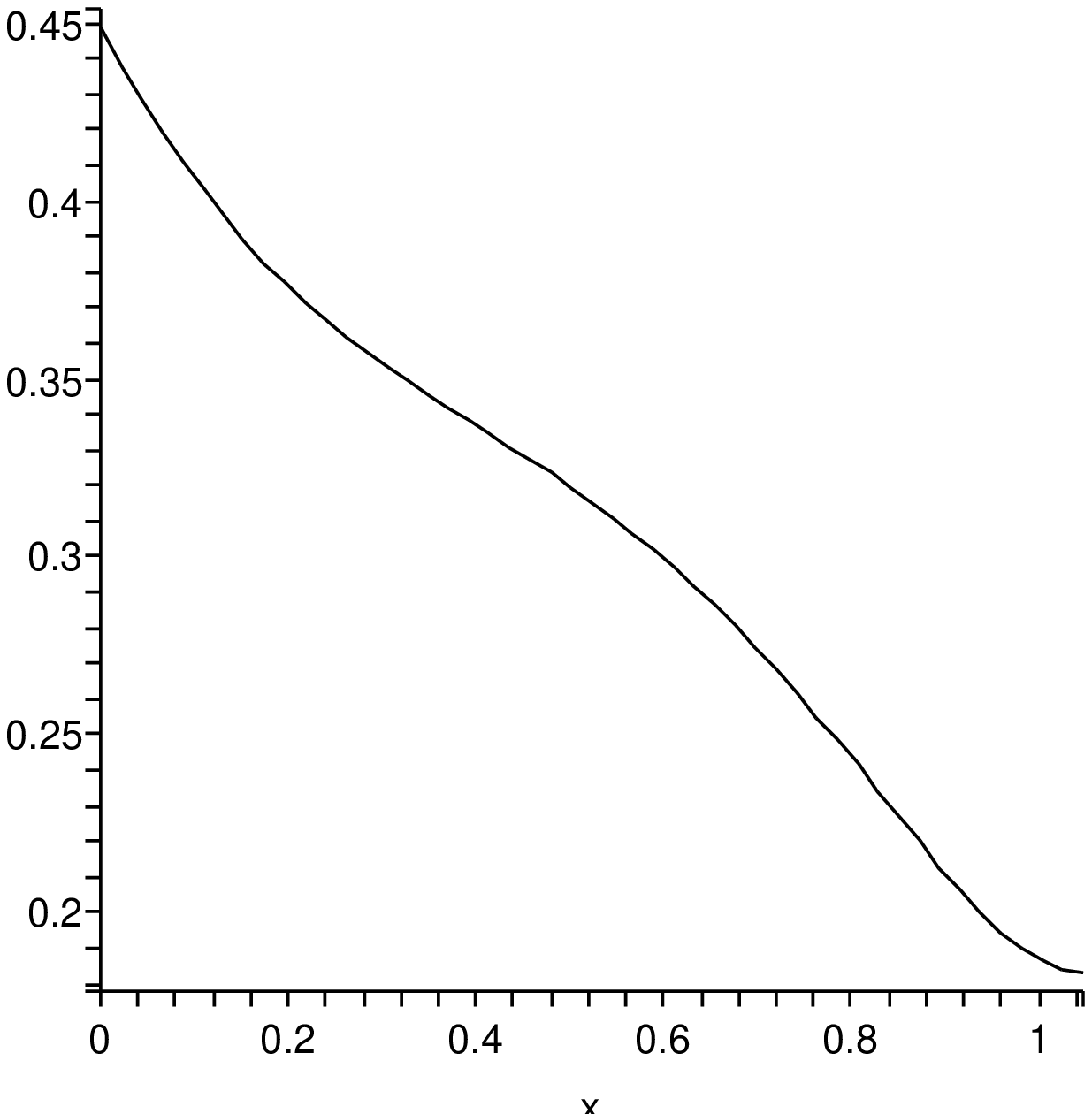}}\hfil
\centerline{Figure 2}
\centerline{Standard Deviation Function}\

The standard deviation function is strictly decreasing (but with an interesting concavity behavior), resulting in the minimum occurring at ${\pi\over3}$, the rightmost endpoint of the interval on which $x$, the length of the arc, is defined.  See Figure 2 above.  The geometric interpretation of the result is that the area of the triangle is $0$, leaving the remaining six regions to each have an area of ${\pi\over6}$.  With six people getting the same size piece of pizza and one person receiving no piece, we refer to this as the ``every seventh child starves" principle.    
Having one child starve is, for us, an unacceptable solution, exposing a limitation of the variance measure in this case.  So let's consider another possible measure - the mean absolute deviation.  The mean deviation function is given by:
\begin{eqnarray*}
{\rm md}(x)&=&{3\over7}\left\vert{x\over2}-2\sin{\lr{x\over2}\rr}\cos{\lr{\pi\over3}+{x\over2}\rr}-{\pi\over7}\right\vert\\
&&{}+{1\over7}\left\vert{-3\sqrt{3}}\cos^2{\lr{\pi\over3}+{x\over2}\rr}+{\pi\over7}\right\vert\\
&&{}+{4\over49}\pi-{3\over14}x+{6\over7}\sin\lr{x\over2}\rr\cos{\lr{\pi\over3}+{x\over2}\rr}-{{3\sqrt{3}}\over7}\cos^2{\lr{\pi\over3}+{x\over2}\rr}.
\end{eqnarray*}

The graph of the absolute deviation function is shown in Figure 3.  The absolute minimum on the interval $[0,{\pi\over3}]$ occurs when the arc length is approximately .96976, resulting in the area of the triangle being a rather small 0.00779, with each circular triangle having an area of 0.44880, and the circular trapezoids each having area 0.59581.  This result is not particularly more satisfying than the fact that minimum SD leads to the ``starving child syndrome."  Looking at the graph there is also a {\it local} minimum.  Exploring this point we find that this value occurs with an arc length of 0.45061 yielding areas for the triangle, circular triangle, and circular trapezoid of 0.44880, 0.09399, and 0.80361, respectively.  But this fact is irrelevant since local minima have no particular significance.

\bigskip
\hfil\scalebox{1}{\includegraphics{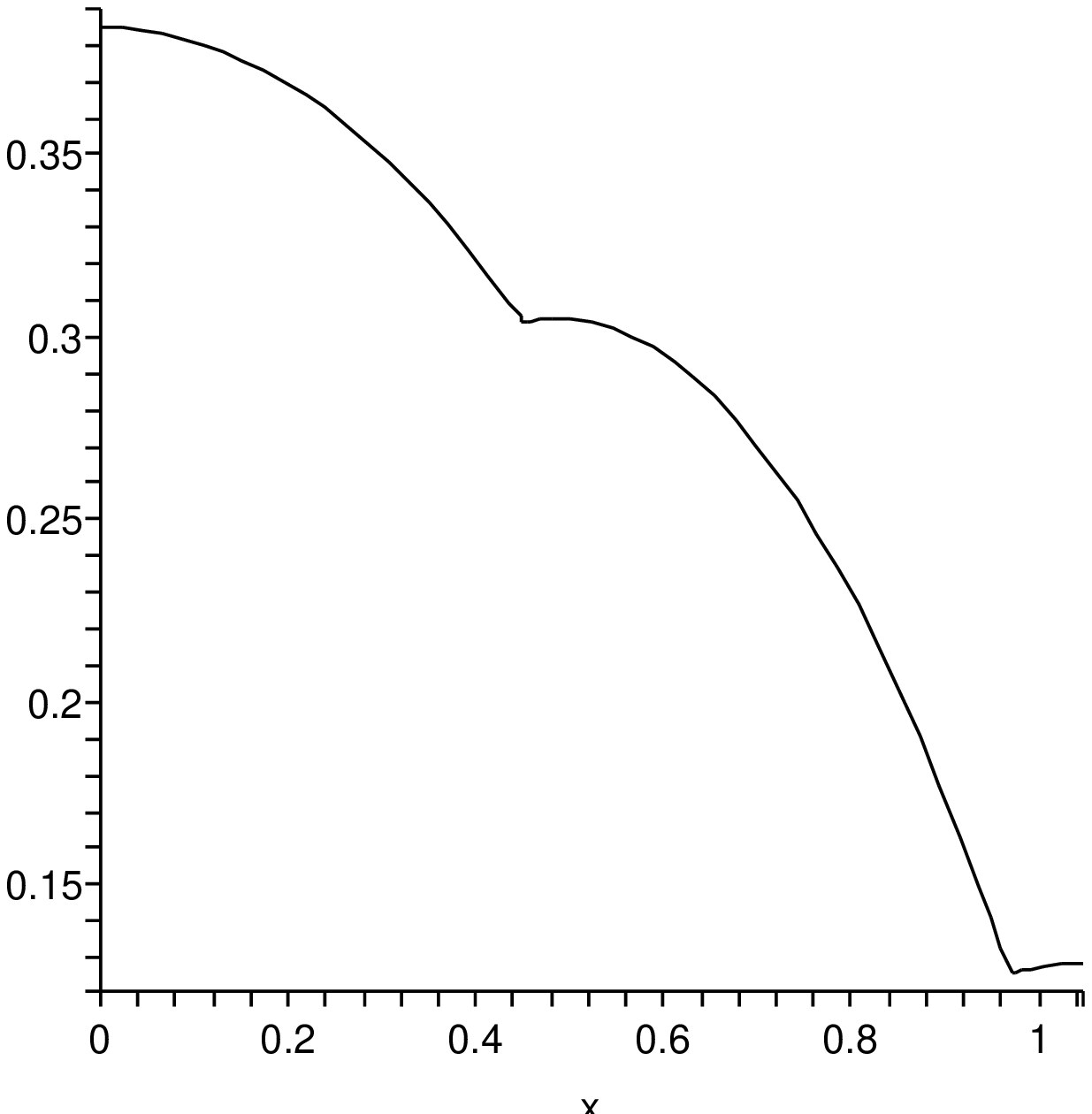}}\hfil
\centerline{Figure 3}
\centerline{Absolute Deviation Function}
\bigskip

%
Still not satisfied, let's explore a third function.  Suppose we seek to maximize the size of the smallest piece, subject to the cuts being of the required kind.  Figure 4 below shows the three area functions graphed on the same set of axes, and it is clear from that diagram that our criterion is met when the triangle and circular triangles each have an area of $\sim 0.2$ and the circular trapezoids are each of area $\sim 0.78$.  This is the solution in which three people will get equal amounts and four people will get equal, but lesser amounts.  We consider this to be our most equitable solution; readers will undoubtedly think of other criteria that might be used in this context.

\bigskip
\hfil\scalebox{1}{\includegraphics{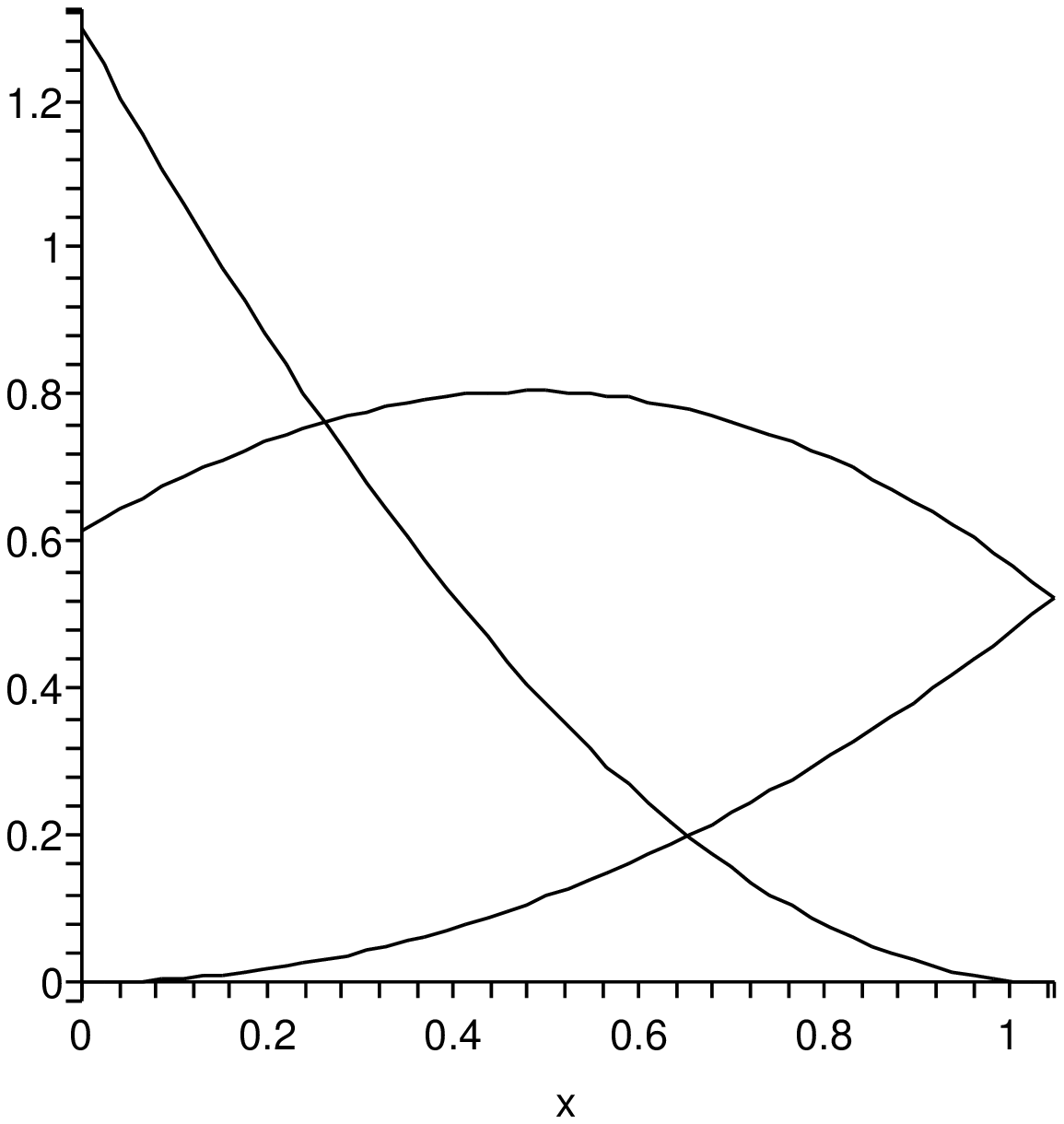}}\hfil
\centerline{Figure 4}
\centerline{The Three Area Functions}
\bigskip

\bigskip

\section{Random Division and Stein's Method}

\nocite{steiner}
\nocite{zimmerman}
\nocite{buck}
\nocite{snell}
\nocite{wetzel}
\nocite{kerr}
\nocite{price}
\nocite{rinott}
\begin{prop}
If $n$ cuts are independently made in a pizza, each with a probability $p$ of success, and so that each cut intersects each prior successful cut, then the expected number of regions is given by
\[\e(R)=1+np+{n\choose 2}p^2.\]
\end{prop}

\noindent{\bf Proof} The number of regions formed by $n$ cuts is given by $1+{n\choose 1}+{n\choose 2}$ or $1+n+{n(n-1) \over 2}$.  Since any given cut will be made with probability $p$ and not made with probability $1-p$, the random variable $X$, denoting the number of cuts made, will possess a binomial distribution.  Thus the number of regions will be
\[R= 1+X+{X(X-1) \over 2},\]
so that 
\begin{eqnarray*}
\e(R)&=&1+\e(X)+\e\lr{X(X-1)} \over 2\rr\nonumber\\
&=&1+np+\sum\lr{{x(x-1)} \over 2}{n\choose x}p^xq^{n-x}\rr\nonumber\\
&=&1+np+{1\over2}p^2{n(n-1)}\sum{(n-2)!\over{(x-2)!(n-x)!}}p^{x-2}q^{n-x}\nonumber\\
&=&1+np+{{n(n-1)}\over2}p^2\nonumber\\
&=&1+np+{n\choose 2}p^2,
\end{eqnarray*}
as claimed.
\begin{lemma}
\begin{eqnarray*}\e(R^2)&=&
\huge(1+3np+{9\over2}n(n-1)p^2+2n(n-1)(n-2)p^3\\&&{}+{1\over4}n(n-1)(n-2)(n-3)p^4\huge)
\end{eqnarray*}
\end{lemma}

\noindent{\bf Proof}  We omit the details.  We use the fact that (i) $\e(Z^2)=\e(Z(Z-1))-\e(Z)$, and (ii) $R-1=X+X(X-1)/2$, as well as formulas for the third and fourth moments of the binomial distribution.
\begin{thm} Assuming that $np\to\infty$ as $n\to\infty$,
$\v(R)=n^3p^3(1-p)(1+o(1))$.
\end{thm}
\noindent{\bf Proof}  
\begin{eqnarray*}
\v(R)&=&\e(R^2)-\lr\e(R)\rr^2\nonumber\\
&=&\bigg(1+3np+{9\over2}n(n-1)p^2+2n(n-1)(n-2)p^3\\
&&{}+{1\over4}n(n-1)(n-2)(n-3)p^4\bigg)\nonumber\\
&{}&-\lr1+2np+n(n-1)p^2+(np)^2+n^2(n-1)p^3+{1\over4}n^2(n-1)^2p^4\rr\nonumber\\
&=&np+{1\over2}n(5n-7)p^2+n(n-1)(n-4)p^3-{1\over2}n(n-1)(2n-3)p^4.\nonumber\\
\end{eqnarray*}
If $np_n=np\to\infty$ as $n\rightarrow\infty$ and if $p$ tends to a non-zero constant, then the last two terms above dominate the variance, and the assertion follows.  If $p\to0$, the third term dominates and the conclusion still holds.

\medskip

Next consider a 3 dimensional object such as a watermelon.  We know that for $n$ cuts there will be a maximum of $1+n+{n\choose2}+{n\choose3}$ pieces \cite{zimmerman}.  As with the pizza, there are $n$ people to make a cut in the watermelon, where each cut will independently be made with probability $p$.  Thus
\[R_{3D}=1+X+{{X(X-1)}\over2}+{{X(X-1)(X-2)}\over6},\ {\rm where}\ X\sim {\rm Bin}(n,p);\]
\begin{eqnarray*}
\e_{3D}(R)&=&1+\e(X)+\e\lr{{{X(X-1)}\over2}}\rr+\e\lr{{{X(X-1)(X-2)}\over6}}\rr\nonumber\\
&=&1+np+{{n(n-1)}\over2}p^2+{1\over6}{n(n-1)(n-2)p^3},
\end{eqnarray*}
after doing a calculation similar to the one above.  Also,
\begin{eqnarray*}
\v_{3D}(R)&=&\e(R^2)-\lr\e(R)\rr^2\nonumber\\
&=&\bigg(1+3np+{9\over2}n(n-1)p^2+{25\over6}n(n-1)(n-2)p^3\nonumber\\
&{}&+{25\over12}n(n-1)(n-2)(n-3)p^4\\&{}&+{5\over12}n(n-1)(n-2)(n-3)(n-4)p^5\nonumber\\
&{}&+{1\over36}n(n-1)(n-2)(n-3)(n-4)(n-5)p^6\bigg)\\
&{}& -\bigg(1+2np+n(2n-1)p^2+
{2\over3}n(n-1)(2n-1)p^3\\&{}{}&+{1\over12}n^2(n-1)(7n-11)p^4
+{1\over6}n^2(n-1)^2(n-2)p^5\\&{}{}&+{1\over36}n^2(n-1)^2(n-2)^2p^6\bigg)\nonumber\\
&=&np+{1\over2}n(5n-7)p^2+{1\over6}n(n-1)(19n-50)p^3\\
&{}&+{1\over2}n(n-1)(3n^2-19n+25)p^4\\
&{}&+{1\over4}n(n-1)(n-2)(n^2-11n+20)p^5\\&{}&-{1\over12}n(n-1)(n-2)(3n^2-15n+20)p^6\nonumber\\
&=&{1\over4}n^5p^5(1-p)(1+o(1))
\end{eqnarray*}
if $np\to\infty$.

Continuing the pizza/watermelon analogy, let's consider a $d$-dimensional ``hypermelon" maximally cut by $n$ hyperplanes.  It has once again been shown by Steiner \cite{steiner} that there would exist $\sum_{i=0}^{d} {n\choose i}$ regions formed \cite{zimmerman}.  Other approaches to the proof of this may be found in \cite{buck}, \cite{kerr}, \cite{price}, and \cite{wetzel}.  If each of the $n$ possible cuts are again successfully made independently with probability $p$,
it may be shown as above that the expected value for the number of regions created in the $d$-dimensional``hypermelon" is given by:
\[\e_{dD}(R)={n\choose 1}p+{n\choose 2}p^2+{n\choose 3}p^3+...+{n\choose d}p^d;\]
of course, in hindsight, the expected value calculation in each dimension follows from the fact that the $r$th factorial moment of a binomial $(n,p)$ random variable equals $n(n-1)\ldots (n-r+1)p^r$.

We now investigate the nature of the concentration of $R$ around its expected value in two dimensions:  By Chebychev's Inequality,
\[P\lr\vert R-\e(R)\vert\ge\lambda\rr \le{{\v(R)}\over{\lambda^2}}\] for every $\lambda>0$, so that we see that 
$P\lr\vert R-\e(R)\vert\ge\lambda\rr\to0$ if $\lambda\gg (np)^{3/2}{\sqrt{1-p}}$, i.e., if $\lambda/ ((np)^{3/2}{\sqrt{1-p}})\to\infty$.  As a result we have, for $p=1/2$, say, that
$R$ is concentrated in an interval of length $\phi(n)n^{3/2}$ around $\e(R)\sim n^2/8$ where $\phi(n)\to\infty$ is arbitrary.  To give two other examples, first, if $p=1/\sqrt n$ then
$\e(R)\sim n/2$ and $\v(R)\sim n^{3/2}$, so the concentration is in a window of magnitude $\phi(n)n^{3/4}$.  Second, if $p=1-(1/{\sqrt n})$, then $\e(R)\sim n^2/2$ and $\v(R)\sim n^{5/2}$, so that the distribution of $R$ is concentrated in an interval of size $\phi(n)n^{5/4}$ around $\e(R)$. 

We next prove a central limit theorem for the distribution of $R$.  Not surprisingly, a CLT does not always hold due to the dependencies among the indicator random variables involved.  However, we are able below to fruitfully use a CLT of Rinott \cite{rinott} to establish the following result: 
\begin{thm} The distribution of $R = R_n$ satisfies, for $p(1-p)^{1/3}\gg n^{-1/9}$,
\[\lim_{n\to\infty}\sup_{w\in{\mathbb R}}\left\vert\p\lr{{R-\e(R)}\over{\sqrt{\v(R)}}}\le w\rr-\Phi(w)\right\vert=0,\]
where $\Phi$ denotes the standard normal distribution function.  (This is tantamount to the assertion that a standardized version of $1+X+X(X-1)/2$ is asymptotically normal when $X\sim{\rm Bi}(n,p)$.)
\end{thm}

\noindent{\bf Proof} Since $X=\sum_{j=1}^nI_j$, where $I_j=1$ iff the $j^{\rm th}$ cut is made and equals zero otherwise, we have
\begin{eqnarray*}
R&=&1+X+{{X(X-1)}\over2}\\
&=&1+\sum_{j=1}^{n} I_j+{1\over2}\lr I_1+...+I_n\rr\lr I_1+...+I_n-1\rr\\
&=&1+{{I_1+...+I_n}\over2}+{1\over2}\lr I_1+...+I_n\rr\lr I_1+...+I_n\rr\\
&=&1+I_1+...I_n+{1\over2}\sum_{i\not=j}I_iI_j
\end{eqnarray*}

We use a bound given by Rinott \cite{rinott} using dependency graphs and Stein's Method. Let the $n^2+1$ summands that constitute $R$ be the {\it vertices} of the graph.  Call these summands $\{J_i: 1\le i\le n^2+1\}$ for convenience.  Draw an {\it edge} between two vertices if they are dependent, i.e. if they share at least one common indicator index $I_r$. Then the {\it degree} of the summand ``1" is zero; for each $j$, the vertex $I_j$ has degree $2(n-1)$; and for any $i,j$, the vertex $I_iI_j$ has degree $4(n-1)$.  Furthermore for any two sets of vertices $A$ and $B$ such that there is no edge beween one vertex in $A$ and another in $B$, the sets $\{J_i:i\in A\}$ and and $\{J_i:i\in B\}$ are independent.  We thus have a dependency graph of $n^2+1$ variables each bounded by $B=1$, and having a {\it maximum degree} of less than $D=4n$.  Rinott's Theorem 1.2 applies to give us, for each $w$,
\begin{eqnarray*}&&{}\left\vert\p\lr{{R-\e(R)}\over{\sigma(R)}}\le w\rr-\int_{-\infty}^{w}{{1}\over{\sqrt{2\pi}}}e^{-{{1}\over2}z^2}dz\right\vert\\
&{}&\le C\max\lc{{ND^2B^3}\over{\sigma^3}}, {{(ND^3B^4)^{1/2}}\over{\sigma^2}},{{DB}\over{\sigma}}\rc,\end{eqnarray*}
where $N$, the total number of summands is bounded by a constant times $n^2$. 
Applying Rinott's Theorem and ignoring constants,
\begin{eqnarray*}
{{ND^2B^3}\over{\sigma^3}}&\sim&{{n^2\cdot n^2}\over{(n^{3}p^{3}(1-p))^{3/2}}}={{1}\over{n^{1/2}p^{9/2}(1-p)^{3/2}}}\to0\\ &{\rm if }&\ p(1-p)^{1/3}\gg{{1}\over{n^{1/9}}};\nonumber\\
{{(ND^3B^4)^{1/2}}\over{\sigma^2}}&\sim&{{(n^2\cdot n^3)^{1/2}}\over{n^3p^3(1-p)}}={1\over{n^{1/2}p^3}(1-p)}\to 0\\ &{\rm if }&\ p(1-p)^{1/3}\gg{1\over {n^{1/6}}};\nonumber\\
{{DB}\over\sigma}&\sim&{{n}\over{n^{3/2}p^{3/2}(1-p)^{1/2}}}={1\over{n^{1/2}p^{3/2}(1-p)^{1/2}}}\to 0\\ &{\rm if}&\ p(1-p)^{1/3}\gg{1\over{n^{1/3}}}.\nonumber
\end{eqnarray*}
The first of the above three conditions is the strongest of the three, proving Theorem 4.  This condition is satisfied if $p$ is fixed or tends to zero so that $p\gg n^{-1/9}$.  If however, $p$ is close to one (and in particular approaches 1 as $n\to\infty$) then the $(1-p)$ term dominates, and we need $(1-p)\gg n^{-1/3}$ for the result to hold.  The cases $p=0,1$ correspond to the deterministic situations in which the CLT is not valid; Theorem 4 gives bounds on how close $p$ may be to either extreme while still giving normal convergence.

\section{(Open) Problems for Research and Classroom Discussion}  (i) How do the results of Section 2 generalize for $n=4,5,\ldots$? (ii) Can we assert that the solution for $n=3$ also minimizes the variance of the pieces among {\it all} possible choices of 3 maximal cuts, both regular (as considered by us) and irregular?
\medskip

\noindent{\bf Acknowledgment}  The research of the both authors was  supported by NSF Grants DMS-9619889 and DMS-0139291. We thank Jeff Knisley for help with our MAPLE calculations.

\end{document}